\documentclass[aos,MSNbibl,seceqn,dvips]{arximspdf}
\usepackage{graphicx}

\doi{10.1214/13-AOS1102} 
\volume{41}
\issue{2}
\pubyear{2013}
\firstpage{870}
\lastpage{896}

\makeatletter

\renewcommand{\mid}{|}
\newcommand{\rrvert}{\vert}
\newcommand{\llvert}{\vert}

\newtheorem{theorem}{Theorem}[section]
\newtheorem{lem}{Lemma}[section]
\newtheorem{prop}{Proposition}[section]

\newproclaim{defn}{Definition}[section]
\newproclaim{rmk}{Remark}[section]
\newproclaim{exmp}{Example}[section]

\makeatother

\begin{document}
\begin{frontmatter}

\title{Bayesian nonparametric analysis of reversible Markov chains}
\runtitle{Bayesian analysis of reversible Markov chains}

\begin{aug}
\author[A]{\fnms{Sergio} \snm{Bacallado}\corref{}\thanksref{t1}\ead[label=e1]{sergiob@stanford.edu}},
\author[B]{\fnms{Stefano} \snm{Favaro}\thanksref{t22}\ead[label=e2]{stefano.favaro@unito.it}}
\and
\author[C]{\fnms{Lorenzo} \snm{Trippa}\ead[label=e3]{ltrippa@jimmy.harvard.edu}}
\runauthor{S. Bacallado, S. Favaro and L. Trippa}
\affiliation{Stanford University, University of Torino and Collegio
Carlo Alberto, Moncalieri, and Harvard School of Public Health and
Dana-Farber~Cancer~Institute}
\address[A]{S. Bacallado\\
Department of Statistics\\
Stanford University \\
Clark Center, S296\\
Stanford, California 94305 \\
USA\\
\printead{e1}}
\address[B]{S. Favaro\\
Department of Economics\\
\quad and Statistics\\
University of Torino\\
Corso Unione Sovietica 218/bis \\
Torino, 10134 \\
Italy\\
\printead{e2}}
\address[C]{L. Trippa\\
Dana-Farber Cancer Institute \\
450 Brookline Ave. CLSB 11039 \\
Boston, Massachusetts 02215 \\
USA\\
\printead{e3}} 
\end{aug}

\thankstext{t1}{Supported by Grant NIH-R01-GM062868 and NSF Grant
DMS-09-00700.}

\thankstext{t22}{Supported by the European Research Council (ERC)
through StG ``N-BNP'' 306406.}

\received{\smonth{2} \syear{2012}}
\revised{\smonth{2} \syear{2013}}

%
\begin{abstract}
We introduce a three-parameter random walk with reinforcement, called
the $(\theta,\alpha,\beta)$ scheme, which generalizes the linearly edge
reinforced random walk to uncountable spaces. The parameter $\beta$
smoothly tunes the $(\theta,\alpha,\beta)$ scheme between this edge
reinforced random walk and the classical exchangeable two-parameter
Hoppe urn scheme, while the parameters $\alpha$ and $\theta$ modulate
how many states are typically visited. Resorting to de Finetti's
theorem for Markov chains, we use the $(\theta,\alpha,\beta)$ scheme to
define a nonparametric prior for Bayesian analysis of reversible Markov
chains. The prior is applied in Bayesian nonparametric inference for
species sampling problems with data generated from a reversible Markov
chain with an unknown transition kernel. As a real example, we analyze
data from molecular dynamics simulations of protein folding.
\end{abstract}

%

\end{frontmatter}

\section{Introduction}
\label{introduction}

The problem that motivated our study is the analysis of benchtop and
computer experiments that produce dynamical data associated with the
structural fluctuations of a protein in water. Frequently, the physical
laws that govern these dynamics are time-reversible. Therefore, a
stochastic model for the experiment should also be reversible.
Reversible Markov models in particular have become widespread in the
field of molecular dynamics~\cite{Pande2010}. Modeling with
reversible Markov chains is also natural in a number of other disciplines.

We consider the setting in which a scientist has a sequence of states
$X_1,\ldots,X_n$ sampled from a reversible Markov chain. We propose a
Bayesian model for a reversible Markov chain driven by an unknown
transition kernel. Problems one can deal with using our model include
(i) predicting how soon the process will return to a specific state of
interest and (ii)
predicting the number of states not yet explored by $X_1,\ldots,X_n$ %
that appear in the next $m$ transitions $X_{n+1},\ldots,X_{n+m}$.
More generally, the model can be used to predict any characteristic of
the future trajectory of the process. Problems (i) and (ii) are of
great interest in the analysis of computer experiments on protein dynamics.

Diaconis and Rolles~\cite{Diaconis2006gd} introduced a conjugate
prior for Bayesian analysis of reversible Markov chains. This prior is
defined via de Finetti's theorem for Markov chains
\cite{Diaconis1980hb}. The predictive distribution is that of a
linearly edge-reinforced random walk (ERRW) on an undirected graph
\cite{Diaconis1988}. Much is known about the asymptotic properties of
this process~\cite{Keane2000df}, its uniqueness~\cite{Rolles2003xd}
and its recurrence on infinite graphs (\cite{Merkl2009}, and references
therein). Fortini, Petrone and Bacallado recently discussed
other examples of Markov chain priors constructed through
representation theorems~\cite{fortini2012hierarchical,Bacallado2011}.

Our construction can be viewed as an extension of the ERRW defined on
an infinite space.
The prediction for the next state visited by the process is not solely
a function of
the number of transitions observed in and out of the last state.
In effect, transition probabilities out of different states share
statistical strength.
This will become relevant in applications where many states are
observed, especially for those states that occur rarely.

A major goal in our application is the prediction of the number of
states that the Markov chain has not yet visited that will appear in
the next $m$ transitions.
More generally, scientists are interested in predicting aspects of the
protein dynamics that may be strongly correlated with the rate of
discovery of unobserved states, for instance, the variability of the
time needed to reach a conformation of interest $y$, starting from a
specific state $x$. Predictive distributions for such attributes are
useful in deciding whether one should continue a costly experiment to
obtain substantial additional information on a kinetic property of interest.

Estimating the probability of discovering new species is a
long-standing problem in statistics~\cite{Bunge1993}. Most
contributions in the literature assume that observations, for example,
species of fish captured in a lake, can be modeled as independent and
identically distributed random variables with an unknown discrete
distribution. In this setting, several Bayesian nonparametric models
have been studied \mbox{\cite{Lijoi2007a,Lijoi2007b,Favaro2009}}. Here we
assume that species, in our case protein conformational states, are
sampled from a reversible Markov chain. To the best of our knowledge,
this is the first Bayesian analysis of species sampling in this setting.

\begin{figure}

\includegraphics{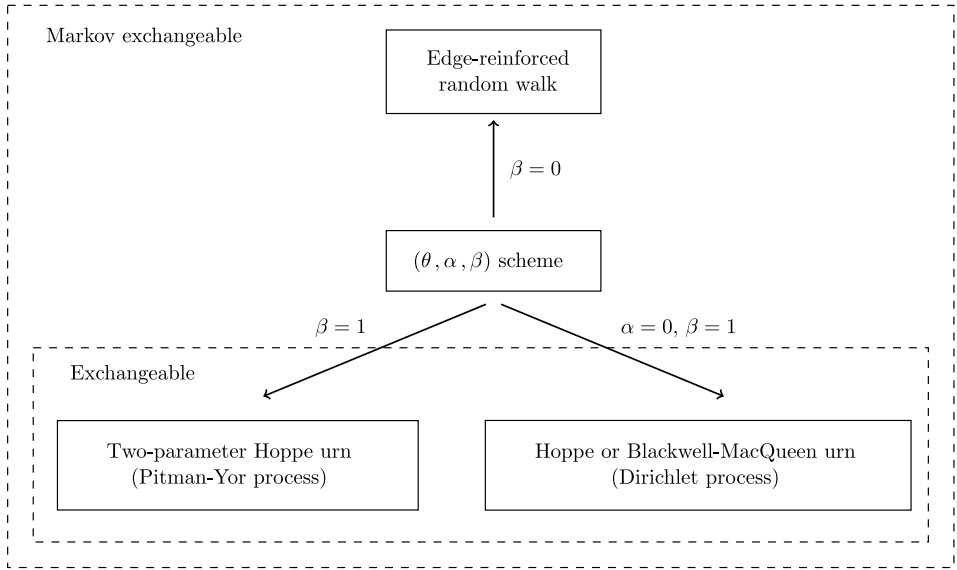}

\caption{Diagram of the $(\theta,\alpha,\beta)$ scheme and special cases.}
\label{diagram}
\end{figure}

We can now outline the article. Section~\ref{tabschemesection}
introduces the species sampling model, which we call the
$(\theta,\alpha,\beta)$ scheme. The process specializes to the ERRW,
a~Markov exchangeable scheme, and to the two-parameter Hoppe urn, a
classical exchangeable scheme which gives rise to the Pitman--Yor
process and the two-parameter Poisson--Dirichlet distribution \cite
{Pitman1996,Pitman1997}. As illustrated in Figure~\ref{diagram}, the
parameter $\beta$ smoothly tunes the model between these two special
cases. Section~\ref{representationssection} shows that the
$(\theta,\alpha,\beta)$ scheme can be represented as a mixture of
reversible Markov chains. This allows us to use its de Finetti measure
as a prior for Bayesian analysis. Section~\ref{largesupportsection}
shows that our scheme is a projection of a conjugate prior for a random
walk on a multigraph. This representation is then used to prove that
our model has full weak support. Section~\ref{sufficientnesssection}
provides a sufficientness characterization of the proposed scheme. This
result is strictly related to the characterizations of the ERRW and the
two-parameter Hoppe urn discussed in~\cite{Rolles2003xd} and
\cite{zabell2005symmetry}, respectively. In Section~\ref{lawsection},
an expression for the law of the $(\theta,\alpha,\beta)$ scheme is
derived, and this result is used in Section
\ref{bayesianinferencesection} to define algorithms for posterior
simulation. Section~\ref{applications} applies our model to the
analysis of two molecular dynamics datasets.
We evaluate the predictive performance of the model by splitting the
data into training and validating datasets. Section \ref
{discussionsection} concludes with a discussion of remaining challenges.


\section{\texorpdfstring{The $(\theta,\alpha,\beta)$ scheme}{The (theta, alpha, beta) scheme}}
\label{tabschemesection}

\begin{figure}

\includegraphics{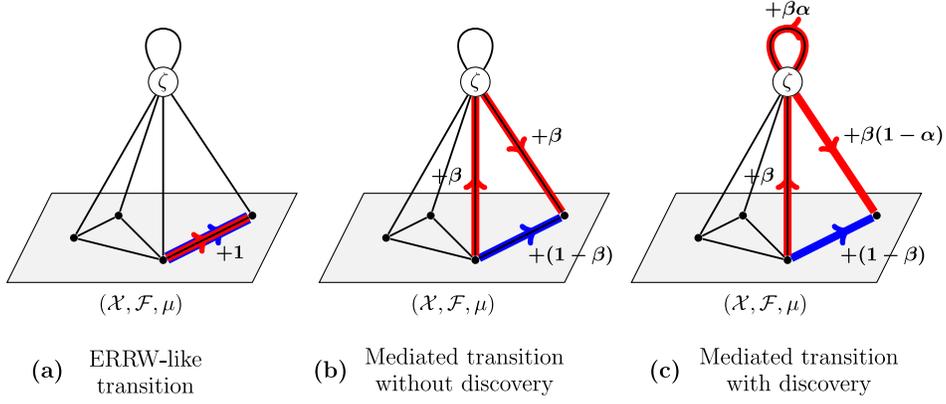}

\caption{Three kinds of transition in the $(\theta,\alpha,\beta)$
scheme. The blue arrow represents the transition between two states in
$\mathcal{X}$, while the red arrows represent the path of an auxiliary
random walk with reinforcement. The edges that have positive weight $g$
before the transition are drawn in black, and in each case, we mark the
reinforcements of $g$ produced by the transition. Self-transitions
follow a slightly different reinforcement scheme formalized in
Definition \protect\ref{tabscheme}.} \label{mechanism}
\end{figure}

The $(\theta,\alpha,\beta)$ scheme is a stochastic process
$(X_i)_{i\in\mathbb{N}}$ on a Polish measurable space $(\mathcal
{X},\mathcal{F})$ equipped with a diffuse (i.e., without point masses)
probability measure $\mu$. We construct the law of the process using an
auxiliary random walk with reinforcement on the extended space
$\mathcal{X}_+:= \mathcal{X}\cup \{\zeta\}$. The auxiliary process
classifies each transition $X_i \to X_{i+1}$ into three categories
listed in Figure~\ref{mechanism} and defines latent variables
$(U_i)_{i\in\mathbb{N}}$, taking values in $\{a,b,c\}$, that capture
each transitions' category. In this section we first provide a formal
definition of the $(\theta,\alpha,\beta)$ scheme and then briefly
describe the latent process.

The law of the $(\theta,\alpha,\beta)$ scheme is specified by a
weighted undirected graph $g$ with vertices in $\mathcal{X}_+$.
This graph can be formalized as a symmetric function $g\dvtx \mathcal
{X}_+\times\mathcal{X}_+\to[0,\infty)$, where $g(x,y)$ is the
weight of an undirected edge with vertices $x$ and $y$. We require that
the set
$\mathcal{S}:= \{ x\in\mathcal{X}; g(\zeta,x)>0\}$
is countable, $\sum_{x\in\mathcal{S}}g(\zeta, x)<\infty$ and that
the set of edges $\{(x,y) \in\mathcal{X}^2; g(x,y)>0 \}$ is a
finite subset of $\mathcal{S}^2$.
The graph will be sequentially reinforced after each transition of the
$(\theta,\alpha,\beta)$ scheme.
In the following definition, we assume the initial state $X_1$ is
deterministic and contained in $\mathcal{S}$.

\begin{defn}
\label{tabscheme} The \textit{$(\theta,\alpha,\beta)$ scheme},
$(X_i)_{i\in\mathbb {N}}$, has parameters $\theta\ge0$,
$\alpha\in[0,1)$ and $\beta\in [0,1]$. The parameter $\theta$ is equal
to the initial weight $g(\zeta,\zeta)$. Suppose we have sampled
$X_1,\ldots,X_i$, where $i\ge1$. Then, given $X_i$ and the reinforced
graph $g$, we sample the following:
\begin{longlist}[(c)]
\item[(a)]
\textit{an ERRW-like transition} to $X_{i+1}$ with probability
\[
\frac{g(X_i,X_{i+1})}{\sum_{x\in\mathcal{X}_+}g(X_i,x)}
\]
and make the following reinforcement:
\[
g(X_i,X_{i+1}) \to g(X_i,X_{i+1}) +1
+\mathbf{1}_{X_i=X_{i+1}};
\]
\item[(b)]
\textit{a mediated transition without discovery} to $X_{i+1}$ with probability
\[
\frac{g(X_i,\zeta)}{\sum_{x\in\mathcal{X}_+}g(X_i,x)} \times\frac
{g(\zeta,X_{i+1}) + \beta\times\mathbf{1}_{X_{i}=X_{i+1}}}{\beta
+\sum_{x\in\mathcal{X}_+}g(\zeta, x)}
\]
and make the following reinforcements:
\begin{eqnarray*}
g(X_i,X_{i+1}) &\to& g(X_i,X_{i+1})
+(1 - \beta)\times(1+\mathbf{1}_{X_i=X_{i+1}}),
\\
g(X_i,\zeta) &\to& g(X_i,\zeta) +\beta,
\\
g(X_{i+1},\zeta) &\to& g(X_{i+1},\zeta) +\beta;
\end{eqnarray*}
\item[(c)]
\textit{a mediated transition with discovery} to a new state $X_{i+1}\sim
\mu$ with probability
\[
\frac{g(X_i,\zeta)}{\sum_{x\in\mathcal{X}_+}g(X_i, x)} \times\frac
{g(\zeta,\zeta)}{\beta+\sum_{x\in\mathcal{X}_+}g(\zeta, x)}
\]
and make the following reinforcements:
\begin{eqnarray*}
g(X_i,X_{i+1}) &\to& g(X_i,X_{i+1})
+(1 - \beta),
\\
g(X_i,\zeta) &\to& g(X_i,\zeta) +\beta,
\\
g(X_{i+1},\zeta)&\to& g(X_{i+1},\zeta) +(1-\alpha)\beta,
\\
g(\zeta,\zeta) &\to& g(\zeta,\zeta) +\alpha\beta.
\end{eqnarray*}
\end{longlist}
\end{defn}

In several applications one may prefer to set the initial $g$ to zero
everywhere except for $g(\zeta,\zeta)=\theta$ and $g(X_1,\zeta)$,
which can be made infinitesimally small. This reduces difficulties
associated with the model specification and does not affect the main
properties of the model discussed in this article. In some cases, we
will relax the assumption that $X_1$ is deterministic, by specifying a
distribution, say $\mu$, for $X_1$ and choosing a positive value for
$g(X_1,\zeta)$. In any case, the conditional distributions $p(X_n|
X_1,\ldots,X_{n-1})$ are dictated by the reinforced scheme in
Definition~\ref{tabscheme}.

We can now describe the latent reinforced process in order to simplify
the interpretation of the $(\theta,\alpha,\beta)$ scheme.
To sample a transition $X_i\to X_{i+1}$, we first take one step in the
auxiliary random walk from $X_i$.
If we land on a state $x\in\mathcal{X}$ [panel~(a) in Figure \ref
{mechanism}], we set $X_{i+1}=x$ and $U_i=a$. If we land on $\zeta$,
we sample another step of the random walk from $\zeta$. Once more, if
we land on some $x'\in\mathcal{X}$ [panel (b) in Figure \ref
{mechanism}], we set $X_{i+1}=x'$ and $U_i=b$. Otherwise [panel (c) in
Figure~\ref{mechanism}], we sample a new state $X_{i+1}$ from $\mu$
and set $U_i=c$.

\begin{rmk}
\label{hoppeurns2} Assume the initial graph $g$ is null everywhere
except for $g(\zeta,\zeta)=\theta$ and $g(\zeta, X_1)=1$. If $\beta=1$,
$\alpha=0$ and $X_1\sim\mu$, then the process $(X_i)_{i\in\mathbb{N}}$
is a Blackwell--MacQueen urn~\cite{Blackwell1973} with base
distribution $\mu$ and concentration parameter $\theta/2$. In different
words, the process is exchangeable and its directing random measure is
the Dirichlet process~\cite{Ferguson1973}.\vadjust{\goodbreak}
\end{rmk}

\begin{rmk}
\label{hoppeurns3}
Under the assumptions in Remark~\ref{hoppeurns2}, by setting $\alpha
>0$ and $g(\zeta, X_1)=1-\alpha$, the process $(X_i)_{i\in\mathbb
{N}}$ is equal to an urn scheme introduced by Engen~\cite{Engen1978},
known in the species sampling context as the two-parameter Hoppe urn;
see Appendix B~\cite{supplement}. This exchangeable process has been
studied extensively by Pitman and Yor~\cite{Pitman1996,Pitman1997}.
Its directing random measure is the Pitman--Yor process \cite
{Ishwaran2001} with base distribution $\mu$, concentration parameter
$(\theta-\alpha)/2$, and discount parameter $\alpha/2$; the sorted
masses of this random measure have the two-parameter Poisson--Dirichlet
distribution.
Note that the discount parameter of the Pitman--Yor process can be
chosen from the unit interval $[0,1)$, while in our construction
$\alpha/2 <0.5$.
\end{rmk}

\begin{rmk}
\label{edgereinforcedrandomwalk}
When $\beta=\theta=0$, edges connected to $\zeta$ cannot be
reinforced, and the $(\theta,\alpha,\beta)$ scheme specializes to
the ERRW on $\mathcal{S}$; see~\cite{Diaconis1988}.
\end{rmk}

\begin{rmk}
Definition~\ref{tabscheme} brings to mind the
two-parameter hierarchical Dirichlet Process hidden Markov model
(HDP-HMM)~\cite{teh2010hierarchical} and its associated species
sampling scheme, the two-parameter Chinese restaurant franchise.
This process has been used for Bayesian modeling of Markov chains on
infinite spaces.
The predictive distribution can be viewed as a $(\theta,\alpha,\beta
)$ scheme in which the underlying infinite graph has directed edges.
However, the $(\theta,\alpha,\beta)$ scheme is not a special case of
this model and has no equivalent hierarchical construction. This
connection is explained in more detail in Appendix B~\cite{supplement}.
\end{rmk}

The influence of each parameter in the $(\theta,\alpha,\beta)$
scheme can be described as follows. The parameter $\beta$ determines
the Markov character of the model; as it approaches 1, the process
becomes exchangeable. The parameter $\theta$ is related to the
concentration parameter of the two-parameter Hoppe urn, which controls
the mode of the number of states visited in a given number of steps.
The parameter $\alpha$ is related to the discount factor in the
two-parameter Hoppe urn, which controls the distribution of frequencies
of different states. It is worth noting that $\beta$ also controls the
number of states visited, which increases markedly as $\beta$ is made larger.

The recurrence of the ERRW on infinite graphs is far from trivial,
especially for locally connected graphs (\cite{Merkl2009}, and references
therein). However, it is not difficult to prove that the
$(\theta,\alpha,\beta)$ scheme a.s. returns infinitely often to all
visited states, the state $\zeta$ is visited infinitely often a.s.
and, if $\theta>0$, the edge $g(\zeta,\zeta)$ is crossed infinitely
often. This notion of recurrence is stated in the next proposition.

\begin{prop}
\label{recurrence}
The $(\theta,\alpha,\beta)$ scheme is recurrent, that is,
the event $\{\sum_{j>i} \mathbf{1}_{X_i=X_j}>0 \}$ has
probability 1 for
every integer $i$.
When $\theta>0$ or when the set $\mathcal{S}$ is infinite, the number
of states visited is infinite almost surely.
\end{prop}


\section{\texorpdfstring{de Finetti representation of the $(\theta,\alpha,\beta)$ scheme}
{de Finetti representation of the (theta, alpha, beta) scheme}}
\label{representationssection}

Diaconis and Freedman defined a special notion of partial
exchangeability to prove a version of de Finetti's theorem for Markov
chains~\cite{Diaconis1980hb}.

\begin{defn}
\label{markovexchangeability} A stochastic process on a countable space
$\mathcal{X}$ is \textit{Markov exchangeable} if the probability of
observing a path $x_1,\ldots,x_n$ is only a function of $x_1$ and the
transition counts $C(x,y):=|\{ x_i=x,x_{i+1}=y; 1\leq i<n\}|$ for all
$x,y\in\mathcal{X}$.
\end{defn}

\begin{theorem}[(Diaconis and Freedman)]
\label{definettistheoremformarkovchains}
A process is Markov exchangeable and returns to every state visited
infinitely often, if and only if it is a mixture of recurrent Markov chains.
\end{theorem}

The $(\theta,\alpha,\beta)$ scheme takes values in an uncountable
space $\mathcal{X}$, which precludes a direct application of Theorem
\ref{definettistheoremformarkovchains}. We will state a more
general notion of Markov exchangeability and use it to prove a de
Finetti style representation for $(X_i)_{i\in\mathbb{N}}$. We use the
notion of \textit{$x$-block} defined in Diaconis and Freedman \cite
{Diaconis1980hb}; given a recurrent trajectory $x_1,x_2,\ldots\,$, the
$i$th $x$-block, for any state $x$ appearing in $x_1,x_2,\ldots\,$, is
the finite subsequence that starts with the $i$th occurrence of $x$ and
ends before the $(i+1)$th occurrence.

The process $(X_i)_{i\in\mathbb{N}}$ visits every state in $\mathcal
{S}$ infinitely often; in addition, it will discover new species in
$\mathcal{X}$ in steps of the third kind in Figure~\ref{mechanism}.
But the new species are sampled independently from $\mu$, which
motivates expressing $(X_i)_{i\in\mathbb{N}}$ as a function of two
independent processes on the same probability space: $(Z_i)_{i\in
\mathbb{N}}$ which represents the sequence where new species are
labeled in order of appearance, and $(T_i)_{i\in\mathbb{N}}$ which
represents the $\mathcal{X}$-valued locations of each species. These
are formally defined in the sequel.

Let $\mathcal{Z}=\mathcal{X}\sqcup\mathbb{N}$, and let $d$ be a
function that maps a sequence in $\mathcal{Z}\supset\mathcal{S}$ to
a sequence in the disjoint union $\mathcal{S}\sqcup\mathbb{N}$. Each
element of the sequence in $\mathcal{S}$ is mapped to itself, and
those not in $\mathcal{S}$ are mapped to the order in which they
appear in the sequence. Hence, the range of $d$ consists of sequences
where every state $j\in\mathbb{N}$ may only appear after all states
$1,2,\ldots,j-1$ have appeared at least once. For example, if $\mathcal
{X}$ is the unit interval and $\mathcal{S}=\{0.1,0.2,0.3\}$, then
\[
d\dvtx  (0.1,7,4,0.3,7,6,4,4) \mapsto(0.1,1,2,0.3,1,3,2,2).
\]
Define $(Z_i)_{i\in\mathbb{N}}:= d((X_i)_{i\in\mathbb{N}})$, and
let $(T_i)_{i\in\mathbb{N}}$ be a sequence of independent random
variables from $\mu$, with $(T_i)_{i\in\mathbb{N}}$ independent from
$(Z_i)_{i\in\mathbb{N}}$. Then,
\[
(X_i)_{i\in\mathbb{N}}\stackrel{\mathrm{d}} {=} (\bar
X_i)_{i\in\mathbb{N}}\qquad \mbox{where we define } \bar X_i:=
\cases{Z_i, &\quad if $Z_i\in\mathcal{S}$,
\cr
T_{Z_i}, &\quad otherwise.}
\]

\begin{prop}
\label{partialexchangeability}
Take two sequences $x_1,\ldots,x_n$ and $x'_1,\ldots,x'_n$ in $\mathcal
{S}\sqcup\mathbb{N}$ that are fixed points of $d$.
Suppose one can map $x_1,\ldots,x_n$ to $x'_1,\ldots,x'_n$ by a
transposition of two blocks in $x_1,\ldots,x_n$ which both begin in\vadjust{\goodbreak}
$x\in\mathcal{S}\sqcup\mathbb{N}$ and end in $y\in\mathcal
{S}\sqcup\mathbb{N}$, followed by an application of the mapping $d$. Then,
\[
p(Z_1=x_1,\ldots,Z_n=x_n) =p
\bigl(Z_1=x'_1,\ldots,Z_n=x'_n
\bigr).
\]
\end{prop}

\begin{exmp}
Assume again $\mathcal{S}=\{0.1,0.2,0.3\}$. If we set
\[
(x_1,\ldots,x_n)=(0.1,1,2,0.3,3,2,4,0.3)
\]
and
\[
\bigl(x'_1,\ldots,x'_n\bigr)=(0.1,1,2,3,0.3,4,2,0.3),
\]
then, by transposing two blocks in $(x_1,\ldots,x_n)$ that start from
$2$ and finish in $0.3$, we obtain the vector $(0.1,1,2,4,0.3,3,2,0.3)$.
Moreover,
\[
d(0.1,1,2,4,0.3,3,2,0.3)=(0.1,1,2,3,0.3,4,2,0.3)=\bigl(x'_1,\ldots,x'_n\bigr).
\]
Proposition~\ref{partialexchangeability} then implies that
the two probabilities, $p(Z_1=x_1,\ldots,Z_n=x_n)$
and $p(Z_1=x'_1,\ldots,Z_n=x'_n)$, are identical.
\end{exmp}

\begin{rmk}
Note that if the process only visits states in $\mathcal{S}$, as is the
case when $\theta=0$, the statement of Proposition \ref
{partialexchangeability} is equivalent to Markov exchangeability; cf.
Proposition 27 in~\cite{Diaconis1980hb}. This fact, together with
Proposition~\ref{recurrence} is enough to show by a straightforward
application of Theorem~\ref{definettistheoremformarkovchains}
that the $(\theta,\alpha,\beta)$ scheme with $\theta=0$ is a
mixture of
recurrent Markov chains on $\mathcal{S}$.
\end{rmk}

Equipped with this notion of Markov exchangeability for the species
sampling sequence $(Z_i)_{i\in\mathbb{N}}$, we show that $(X_i)_{i\in
\mathbb{N}}$ can be represented as a mixture of Markov chains.

\begin{prop}
\label{representationtheorem}
There exists a mixture of Markov chains $(W_i)_{i\in\mathbb{N}}$,
taking values in $\mathcal{S}\sqcup\mathbb{N}$, such that, if we define
\[
\tilde X_i: = \cases{W_i, &\quad if $W_i\in
\mathcal{S}$,
\cr
T_{W_i}, &\quad otherwise,}
\]
where $(W_i)_{i\in\mathbb{N}}$ and $ (T_i)_{i\in\mathbb{N}}\sim\mu
^{\mathbb{N}}$ are independent, then $(\tilde X_i)_{i\in\mathbb
{N}}\stackrel{\mathrm{d}}{=}(X_i)_{i\in\mathbb{N}}$. That is, for
some measure $\phi$ on $(\mathcal{S}\sqcup\mathbb{N}) \times
\mathcal{P}$, where $\mathcal{P}$ is the space of stochastic matrices
on $\mathcal{S}\sqcup\mathbb{N}$, the distribution of $(W_i)_{i\in
\mathbb{N}}$, can be represented as
\[
p(W_1=w_1,\ldots,W_n=w_n) =
\int_{\mathcal{P}} \prod_{i=1}^{n-1}P(w_i,w_{i+1})
\phi(w_1,dP).
\]
\end{prop}

Let $\mathcal{P}_r\subset\mathcal{P}$ be the set of transition
probability matrices for recurrent reversible Markov chains. We can
show that this set has probability 1 under the de Finetti measure.

\begin{prop}
\label{reversibility}
$\phi((\mathcal{S}\sqcup\mathbb{N}) \times\mathcal{P}_r) = 1$.
\end{prop}


\section{\texorpdfstring{The $(\theta,\alpha,\beta)$ scheme with colors}{The (theta, alpha, beta) scheme with colors}}
\label{largesupportsection}

This section shows that the $(\theta,\alpha,\beta)$ scheme can be
interpreted as a Bayesian conjugate model for a random walk on a
multigraph. This representation is used for showing that the
$(\theta,\alpha,\beta)$ scheme has large support in a sense that will
be made precise in the sequel. We also make a connection between the de
Finetti measure of the ERRW on a finite graph and our model.

We start by defining a colored random walk on a weighted multigraph
$\tilde g$. The vertices of the graph take values in $\mathcal{X}$,
and we now allow there to be more than one edge between every pair of
vertices. Every edge is associated to a distinct color in a set
$\mathcal{C}$. We assign a weight $\tilde g (\{x,y\}, c)$ to the edge
connecting $x$ and $y$ with color~$c$, requiring that $\tilde g(x):=
\sum_{y,c} g(\{x,y\},c) <\infty$ for all $x\in\mathcal{X}$. A
random walk on this graph is a process that starts from $x_1\in
\mathcal{X}$, and after arriving at some state $x$, traverses the edge
$(\{x,z\},c)$ with probability $\tilde g (\{x,z\}, c)/ \tilde g(x)$.
Let $p_{\tilde g}$ be the law of this process.

A Bayesian statistician observes a finite sequence of traversed colored
edges and wants to predict the future trajectory of the colored random
walk. We suggest how to use the $(\theta,\alpha,\beta)$ scheme in
this context. Informally, in a ERRW-like transition, we reinforce a
single edge of a specific color, while in a mediated transition, we
draw a new edge with a novel color.

%

The Bayesian model is a random sequence of colored edges $(E_i)_{i\in
\mathbb{N}}$.
We use $C_i$ and $\{X_i,X_{i+1}\}$ to denote the color and vertices of $E_i$.
Let $\mu$ and $\gamma$ be nonatomic distributions over $\mathcal
{X}$ and $\mathcal{C}$, respectively, and specify $\theta$, $\alpha$
and $\beta$ as in the previous sections.
Let $X_1=x_1$, $C_1\sim\gamma$ and
\[
X_2\mid X_1,C_1 \sim\frac{\beta(1-\alpha) }{\beta+\theta}
\bolds{\delta}_{X_1} + \frac{\theta+\alpha\beta}{\beta+\theta} \mu.
\]
The distribution of $(X_1,X_2)$ corresponds to an initial graph with
only 2 edges, with endpoints $\{x_1,\zeta\}$ and $\{\zeta,\zeta\}$,
weighted by $-\alpha\beta$ and $\theta+\alpha\beta$, respectively.
If the initial weighted graph $g$ in a $(\theta,\alpha,\beta)$
scheme is chosen as above, the reinforcement rules in Definition \ref
{tabscheme} produce a well-defined process even if the initial value
of $g(x_1,\zeta)$ is negative. After the first transition all edges
will have nonnegative weights.
This choice for the initial weighted graph will be used in the present
section and Section~\ref{sufficientnesssection}.

After a path $(E_1,\ldots,E_n)$,
the probability of recrossing an edge $E_j$, with $j\le n$ and $
X_{n+1}\in\{X_{j},X_{j+1}\} $, is
\[
p(E_{n+1} = E_j \mid E_1,\ldots,E_n )=
\frac{G_{n}(\{X_j,X_{j+1}\},C_j)}{W_{X_{n+1},n}},
\]
where
\begin{eqnarray*}
G_{n}\bigl(\{x,y\},c\bigr)&=&\max\biggl(0, -\beta+\sum
_{i\le n} \mathbf{1}_{
E_i=(\{x,y\},c)} \biggr)2^{\mathbf{1}_{x=y}},
\\
W_{x,n}&=&-\alpha\beta+\sum_{i\le n} (
\mathbf{1}_{x\in\{
X_i,X_{i+1}\}} ) 2^{\mathbf{1}_{X_i=X_{i+1}}}.
\end{eqnarray*}
In words, the probability of recrossing an edge is linear in the number
of crossings. Let $\mathcal{C}_n$ be the set of distinct colors in
$(C_1,\ldots,C_n)$. The conditional probability that
$C_{n+1}\notin\mathcal{C}_n$,
\[
p(C_{n+1}\notin\mathcal{C}_n\mid E_1,\ldots,E_n) = \frac
{B_{X_{n+1},n}}{W_{X_{n+1},n}},
\]
where
\[
B_{x,n}=- \alpha\beta+\beta\sum_{i\le n} (
\mathbf{1}_{x\in
\{X_i,X_{i+1}\}} ) ( \mathbf{1}_{C_i \notin\mathcal
{C}_{i-1}} )2^{\mathbf{1}_{X_i=X_{i+1}}}
\]
is linear in the number of distinct colored edges adjacent to $X_{n+1}$.
The probability that $X_{n+2}=y$, for any
vertex $y\in\{X_1,\ldots,X_n\}$, conditional on $C_{n+1}\notin
\mathcal{C}_n$, is
%
\begin{equation}
\label{newedgetooldvertex} p(X_{n+2} = y \mid E_1,\ldots,E_n,
C_{n+1}\notin\mathcal{C}_n) = \frac{ \beta\mathbf{1}_{X_{n+1}=y}+
B_{y,n} }{2\beta\times |\mathcal{C}_n|+\theta+\beta},
\end{equation}
which depends linearly on the number of distinct colored edges adjacent
to~$y$. Finally,
\[
X_{n+2}\mid C_{n+1}\notin\mathcal{C}_n,\qquad X_{n+2}
\notin\{X_1,\ldots,X_{n+1}\},\qquad E_1,\ldots,E_n \sim\mu
\]
and
\[
C_{n+1}\mid X_{n+2},\qquad E_1,\ldots,E_n,\qquad C_{n+1}
\notin\mathcal{C}_n \sim\gamma.
\]
The following property is a direct consequence of this definition.

\begin{prop}\label{eqd}
The sequence $(X_i)_{i\in\mathbb{N}}$ of $\mathcal{X}$-valued states
visited by $(E_i)_{i\in\mathbb{N}}$ is identical in distribution to a
$(\theta,\alpha,\beta)$ scheme initiated at $x_1$, with $g$ everywhere
null except at $g(x_1,\zeta)=-\beta\alpha$ and
$g(\zeta,\zeta)=\theta+\alpha\beta$.
\end{prop}

Furthermore, given a $(\theta,\alpha,\beta)$ scheme with an
arbitrary initial graph $g$, it is possible to construct a colored
$(\theta,\alpha,\beta)$ scheme with a closed-form predictive
distribution such that the equality stated in Proposition~\ref{eqd} holds.
This would require changing the definition of $(E_i)_{i\in\mathbb
{N}}$ in a way that preserves the reinforcement scheme.


\begin{prop}\label{exc}
Let $e_1,\ldots,e_n$ be a colored path, with $e_i=(\{x_i,x_{i+1}\},\allowbreak c_i)$,
and let $\Lambda(e_1,\ldots,e_n)$ be the probability of
the event
\[
\bigcap_{ e_i=e_j} \{E_i=E_j\}
\bigcap_{e_i \neq e_j} \{E_i\neq E_j\}
\bigcap_{x_i=x_j} \{X_i=X_j\}
\bigcap_{x_i\neq x_j} \{X_i\neq X_j
\}.
\]
Suppose $e_{\sigma(1)},\ldots,e_{\sigma(n)}$, for some permutation
$\sigma$, is also a colored path starting at $x_1$. Then,
%
\begin{equation}
\label{excequation} \Lambda(e_1,\ldots,e_n)=
\Lambda(e_{\sigma_1},\ldots,e_{\sigma
_n}).
\end{equation}
\end{prop}

This result, related to Proposition~\ref{markovexchangeability},
establishes a probabilistic symmetry between paths that can be mapped
to each other\vadjust{\goodbreak} by permuting the order of edges crossed, and applying
certain automorphisms to $\mathcal{X}$ and $\mathcal{C}$. Proposition
\ref{exc} gives rise to the following de Finetti representation.

\begin{prop}
\label{deFinetticolor} There exists a mixture of colored random walks
on weighted multigraphs $(\{W_i,W_{i+1}\},J_i)_{i\in\mathbb{N}}$, with
$W_i\in\mathbb{N}$ and $J_i\in\mathbb{N}$, and independent processes
$(T_i)_{i\in\mathbb{N}}\sim\mu^{\mathbb {N}}$ and
$(V_i)_{i\in\mathbb{N}}\sim\gamma^{\mathbb{N}}$, such that the sequence
of colored edges $(\{\tilde X_i,\tilde X_{i+1}\},\tilde
C_i)_{i\in\mathbb{N}}$, defined by
\[
(\tilde X_i, \tilde C_i) = \cases{ (x_1,
V_{J_i}), &\quad if $W_i = 1$,
\cr
(T_{W_i},
V_{J_i}), &\quad if $W_i >1$,}
\]
is identical in distribution to $(E_i)_{i\in\mathbb{N}}$.
\end{prop}

The proof of Proposition~\ref{deFinetticolor} constructs the
discrete process $(W_i,J_i)_{i\in\mathbb{N}}$, which will now be used
to show that the $(\theta,\alpha,\beta)$ scheme has large support.
The law of $(W_i,J_i)_{i\in\mathbb{N}}$ can be mapped bijectively to
the exchangeable law of the sequence of $x_1$-blocks in the process,
which in the present section and the next are defined as sequences of
labeled edges.
Therefore, the de Finetti measure of $(W_i,J_i)_{i\in\mathbb{N}}$
uniquely identifies the de Finetti measure of the $x_1$-blocks, and
vice versa.
We denote the random distribution of the $x_1$-blocks $\eta$. Note
that the space of $x_1$-blocks is discrete because
$W_i$ and $J_i$ are integer-valued.

We show that the law of $\eta$ has full weak support.
On the basis of Proposition~\ref{eqd}, we then conclude that the
$x_1$-block de Finetti measure induced by the $(\theta,\alpha,\beta
)$ scheme also has full weak support.
The next proposition is proven for $(E_i)_{i\in\mathbb{N}}$ as
defined in this section; the result can be extended to any analogous
reinforced process corresponding to a specific $(\theta,\alpha,\beta
)$ scheme.

\begin{prop}
\label{largesupportprop}
Let $\eta^o$ be the $x_1$-block distribution induced by
a colored random walk on an arbitrary weighted multigraph
with vertices and colors in $\mathbb{N}$, in which the sum of weights
is finite, $\sum_{x\in\mathbb{N}} \tilde g(x)<\infty$.
For every $\varepsilon>0$, $m\ge1$ and any collection of bounded real
functions $f_1,\ldots,f_m$ on the space of $x_1$-blocks,
\[
p \bigl( \eta\in U_{\varepsilon,f_1,\ldots,f_m }\bigl(\eta^o\bigr) \bigr
)>0,
\]
where
\[
U_{\varepsilon,f_1,\ldots,f_m }\bigl(\eta^o\bigr)= \biggl\{\eta^\prime\dvtx
\biggl\llvert\int f_i \,d\eta^o-\int f_i \,d
\eta^\prime\biggr\rrvert<\varepsilon, i=1,\ldots,m \biggr\}.
\]
\end{prop}

The process $(E_i)_{i\in\mathbb{N}}$ also reveals a connection between
the $(\theta,\alpha,\beta)$ scheme and the ERRW.
Recall that the colored $x_1$-blocks $(H_i)_{i\in\mathbb{N}}$ in the
Bayesian model are exchangeable and, by de Finetti's theorem,
conditionally independent.
Consider their posterior distribution given the subsequence
$(E_i)_{i\le n}$, and assume it has $X_{n+1}=X_1$ and includes\vadjust{\goodbreak} $k$
colored $x_1$-blocks.
These assumptions are only made to simplify the exposition.
The limits
%
\begin{equation}
\label{natk} T_{(\{x,y\},c)} = \lim_{m\to\infty} \biggl( \frac{
\sum_{i=n}^{n+m}\mathbf{1}_{(\{X_i,X_{i+1}\},C_i)= ( \{x,y\},c)} }{
\sqrt{ \sum_{i=n}^{n+m}\mathbf{1}_{ X_i=x}} \sqrt{
\sum_{i=n}^{n+m}\mathbf{1}_{ X_i=y}} } \biggr)
\end{equation}
for every edge $(\{x,y\},c)\in\{E_1,\ldots,E_n\}$ are functions of
the directing random measure for the sequence of $x_1$-blocks. From
these limits, one can obtain the probability in the directing random
measure of any $x_1$-block formed with edges in $\{E_1,\ldots,E_n\}$.
Namely, given $E_1,\ldots,E_n$ and the tail $\sigma$-field of
$(E_i)_{i\in\mathbb{N}}$, the probability of an $x_1$-block
$e_1,\ldots,e_k$ is the product $\prod_{i=1}^k T_{e_i}$.
We can now state the connection with the ERRW.

\begin{prop}\label{ide}
There exists an ERRW on a multigraph with a finite number of edges,
such that the joint posterior distribution of the random variables
in~(\ref{natk}), given $E_1,\ldots,E_n$, is identical to the
distribution of the same limits in the ERRW.
\end{prop}
Appendix C~\cite{supplement} contains a constructive proof of this
proposition, in which one such ERRW, whose parameters depend on
$E_1,\ldots,E_n$, is defined.


\section{Sufficientness characterization}
\label{sufficientnesssection}

This section provides a characterization of the colored
$(\theta,\alpha,\beta)$ scheme in terms of certain predictive
sufficiencies or \textit{sufficientness} conditions. The first
characterization of this type, for the P\'olya urn, was proven in an
influential paper by Johnson~\cite{Zabell1982ir} and has since
been extended to other predictive schemes for discrete sequences such
as the two-parameter Hoppe urn~\cite{zabell2005symmetry} and the
edge-reinforced random walk~\cite{Rolles2003xd}. Our result is closely
related to the work of these authors and uses similar proof techniques.
In addition to its clear subjective motivation, the characterization
elucidates connections between the $(\theta,\alpha,\beta)$ scheme and
other popular nonparametric Markov models
\cite{fortini2012hierarchical}.

Consider a random sequence of colored edges $(\tilde E_i)_{i\in\mathbb
{N}}$, where $\tilde E_i=(\{X_i,\break X_{i-1}\},C_i)$.
Each color in $(\tilde E_i)_{i\in\mathbb{N}}$ identifies an edge.
We assume $(\tilde E_i)_{i\in\mathbb{N}}$ is a mixture of random
walks on weighted and colored multigraphs, in the sense of Proposition
\ref{deFinetticolor}, which visits more than 2 vertices with
probability 1.
It will be shown that, if the predictive distribution of $(\tilde
E_i)_{i\in\mathbb{N}}$ satisfies certain conditions, then the process
is a $(\theta,\alpha,\beta)$ scheme.\vspace*{1pt}

Consider a path $Z_n=(\tilde E_1,\ldots,\tilde E_n)$ and define:
\begin{eqnarray*}
\mbox{(i)}\hspace*{10.40pt}\quad  \kappa(Z_n) &=& \bigl|\{0\leq i\leq n; X_i =
X_n\}\bigr| + \frac{\mathbf{1}_{X_n\neq X_0}}{2},
\\
\mbox{(ii)}\quad \kappa(e,Z_n) &=& \bigl|\{1\leq i\leq n; \tilde
E_i = e\} \bigr| \quad\mbox{and}
\\
\mbox{(iii)}\hspace*{10.75pt}\quad \tau(Z_n)&=& \sum_{v}
\eta(v,Z_n)/2
\\
\mbox{(iv)\hspace*{37.5pt}}\quad &&\hspace*{-16pt}
\mbox{where }\eta(v,Z_n) = \sum_{i\le n} (
\mathbf{1}_{v\in\{X_i,X_{i+1}\}} ) ( \mathbf{1}_{C_i
\notin\{C_1,\ldots, C_{i-1}\}} )2^{\mathbf{1}_{X_i=X_{i+1}}}.
\end{eqnarray*}
These variables describe: (i) how many times $X_n$ has been visited and
whether it coincides with $X_0$,
(ii) how many times an edge $e$ has been traversed, (iii) the number of
observed colors and (iv) the degree of $v\in\mathcal{X}$ in
the multigraph constructed by all distinct colored edges in $Z_n$.
In summary, they are easily interpretable. We also use
$ \rho(Z_n) $ to denote the number of distinct $\mathcal{X}$-valued
states in $Z_n $, and
the indicator $s(e)$, which is
equal to 1 if $e$ is a loop and $0$ otherwise.

We can now define sufficientness conditions for $(\tilde E_i)_{i\in\mathbb
{N}}$. The process satisfies Condition 1 if there exist functions $b_0$
and $b_1$ such that for every $ e\in\{ \tilde E_{1},\ldots, \tilde
E_{n}\}$ incident on $X_n$,
%
\begin{equation}
\label{a} p (\tilde E_{i+1}=e | Z_n ) = b_{s(e)}
\bigl(\kappa(Z_n),\kappa(e,Z_n) \bigr) \in(0,1).
\end{equation}
In words, the probability of making a transition through an edge $e=(\{
X_n,v\},\break c)$
in $Z_n$ depends on the number of times the edge has been crossed and
the number of visits to $X_n$. The process satisfies Condition 2 if
there is a function $g$ such that
%
\begin{equation}
\label{b} p \bigl( \tilde E_{n+1}\notin\{ \tilde E_{1},\ldots, \tilde E_{n}\} | Z_n \bigr) = g \bigl(
\kappa(Z_n),\eta(X_n,Z_n) \bigr) \in(0,1).
\end{equation}
That is, the probability of a transition through a new edge is a
function of the number of observed edges $\eta(X_n,Z_n)$ incident on
$X_n$, and the number of visits to $X_n$. Condition 3 requires that
some function $h$ satisfies, for every $v\in\{X_0,\ldots,X_n\}$,
%
\begin{eqnarray}
\label{c}
&& p \bigl( X_{i+1} = v | Z_n, \tilde
E_{n+1}\notin\{ \tilde E_{1},\ldots, \tilde E_{n}
\} \bigr)
\nonumber\\[-8pt]\\[-8pt]
&&\qquad= h \bigl(\tau(Z_n),\eta(v,Z_n) +
\mathbf{1}_{X_n=v} \bigr)\in(0,1).\nonumber
\end{eqnarray}
If a new edge will be traversed, then the conditional probability that
the path will go to an already seen vertex depends solely on the number
of edges out of said vertex and the overall number of observed edges.
Finally, the process satisfies Condition~4 if there is a function $q$,
such that
%
\begin{eqnarray}
\label{d}
&& p \bigl( X_{i+1} \notin\{X_0,\ldots,X_n\} | Z_n, \tilde E_{n+1}\notin\{ \tilde
E_{1},\ldots, \tilde E_{n}\} \bigr)
\nonumber\\[-8pt]\\[-8pt]
&&\qquad =q \bigl(\tau(Z_n),\rho(Z_n) \bigr)\in(0,1);\nonumber
\end{eqnarray}
that is, the conditional probability that the path will go to an
unseen vertex is a function of the total number of edges and vertices.

The main result of this section can be divided into two lemmas, the
first of which depends only on 3 of the conditions above.

\begin{lem}
\label{suff1} If the process $(\tilde E_i)_{i\in\mathbb{N}}$ satisfies
Conditions 1, 2 and 3, there exist $\beta\in[0,1)$ and
$\lambda\in[-\beta,\infty)$ such that
%
\begin{equation}\label{nicerformg}
b_s(k,j) = \frac{(1+s)(j -\beta) }{\lambda+ 2k-2} \quad\mbox{and}\quad g(k,t) =
\frac{\lambda+ \beta t }{ \lambda+ 2k-2}.\vadjust{\goodbreak}  
\end{equation}
\end{lem}

\begin{lem}
\label{suff2}
If the process $(\tilde E_i)_{i\in\mathbb{N}}$ satisfies Conditions
1, 2, 3 and 4, there exist $\alpha\in[0,1)$ and $\lambda'\in
[-\alpha,\infty)$ such that
%
\begin{equation}\label{nicerformg}
h(n,j) = \frac{j -\alpha}{\lambda' + 2n+1-\alpha} \quad\mbox{and}\quad q(n,t) =
\frac{\lambda' + \alpha(t-1) }{ \lambda' + 2n+1-\alpha}
\end{equation}
for $n=1,2,\ldots\,$, $j=1,\ldots,2n$ and $t=2,\ldots,n+1$.
\end{lem}

The characterization of the $(\theta,\alpha,\beta)$ scheme with
colors follows from the two previous lemmas.

\begin{theorem}
\label{sufficientness}
If the process $(\tilde E_i)_{i\in\mathbb{N}}$ satisfies Conditions
1, 2, 3 and 4, then there exist $\alpha\in[0,1)$, $\beta\in[0,1)$
and $\theta>-2\alpha\beta$, such that the conditional transition
probabilities of the process $(\tilde E_i)_{i\in\mathbb{N}}$, given
any path that visits more than 2 vertices, are equal to those in the
$(\theta,\alpha,\beta)$ scheme with colors $(E_i)_{i\in\mathbb{N}}$.
\end{theorem}


\section{\texorpdfstring{The law of the $(\theta,\alpha,\beta)$ scheme}
{The law of the (theta, alpha, beta) scheme}}
\label{lawsection}

We provide an expression for the law of the species sampling sequence
$(Z_i)_{i\in\mathbb{N}}$, defined in Section \ref
{representationssection}. Recall that the process takes values on
$\mathcal{S}\sqcup
\mathbb{N}$, and consider a fixed path $\mathbf{z}$.

Let $n_{xy}$ be the number of transitions in $\mathbf{z}$ between $x$
and $y$ in either direction and $n_x = \sum_{y\in\mathcal{S}\sqcup
\mathbb{N}} n_{xy}$. For each pair $x,y\in\mathcal{S}\sqcup\mathbb
{N}$, we introduce $k_{xy}\leq n_{xy}$ for the number of ERRW-like
transitions [panel (a), Figure~\ref{mechanism}] between $x$ and $y$ in
either direction. Let $\ell_x:= \sum_{y\in\mathcal{S}\sqcup
\mathbb{N}} (n_{xy}-k_{xy})2^{\mathbf{1}_{x=y}}$ be the number of
times that the latent path traverses the edge $(\zeta,x)$ in either
direction, let $\ell:= \sum_{x\in\mathcal{S}\sqcup\mathbb{N}}
\ell_x/2$ be the number of mediated transitions, and let $\ell'$ be
the number of times that $(\zeta,\zeta)$ is traversed. Note that
$\ell$, $\ell_x$ and $\ell'$ are functions of $\mathbf{z}$ and
$\mathbf{k}= \{k_{xy};x,y\in\mathcal{S}\sqcup\mathbb{N}\}$. We
will also need $g(x) = \sum_{y\in\mathcal{S}}g(x,y)$, where $g$ is
the initial weighted graph.

Given $\mathbf{k}$ and $\mathbf{z}$, we know the number of
transitions out of $\zeta$ and out of $x\in\mathcal{S}\sqcup\mathbb
{N}$ in the latent path. Each transition adds a factor to the
denominator of the probability of a latent path, which increase by a
fixed amount, $2\beta$ or~$2$, between occurrences. Similarly, given
$\mathbf{k}$, we know the number of times that $(\zeta,x)$ is
traversed; each transition adds a factor in the numerator, and these
factors are sequentially reinforced by a fixed amount $\beta$.
Finally, $(\zeta,\zeta)$ is traversed $\ell'$ times, and this
contributes a factor $\theta(\theta+\alpha\beta)\cdots(\theta
+[\ell'-1]\alpha\beta)$ to the numerator of the probability of the
latent path.

We can write $p(\mathbf{z}) = \sum_{\mathbf{k}} p(\mathbf
{z},\mathbf{k})$, where $p(\mathbf{z},\mathbf{k})$ is the total
probability of all latent paths consistent with $(\mathbf{z},\mathbf
{k})$. Taking into account that the factors listed in the previous
paragraph are common to all latent paths with a given $(\mathbf
{z},\mathbf{k})$, we obtain
\begin{eqnarray*}
p(\mathbf{z},\mathbf{k})
&=&
F(\mathbf{z},\mathbf{k}) (\theta)_{\ell' \uparrow
\alpha\beta} \prod_{x: n_x>0}
\bigl(g(x,\zeta)+\beta(1-\alpha)\mathbf{1}_{x\in\mathbb{N}}
\bigr)_{\ell_x-\mathbf{1}_{x\in\mathbb{N}} \uparrow\beta} \\
&&\hspace*{0pt}{}\Big/\biggl(
\bigl(g(\zeta)+\beta\bigr)_{{\ell}\uparrow2\beta}
\bigl(g(z_1) \bigr)_{\lfloor(n_{z_1}+1)/2 \rfloor\uparrow2}
\\
&&\hspace*{13.3pt}{}\times\mathop{\prod_{x:n_x>0}}_{x\neq z_1}
\bigl(g(x)+1-\alpha\beta\mathbf{1}_{x\in\mathbb{N}}
\bigr)_{\lfloor
n_{x}/2 \rfloor\uparrow2} \biggr),
\end{eqnarray*}
where we use Pitman's notation for factorial powers
\[
(r)_{n\uparrow q}:= r(r+q) (r+2q)\cdots\bigl(r+(n-1)q\bigr).
\]
The function $F(\mathbf{z},\mathbf{k})$ is a sum with as many terms
as the possible latent paths consistent with $(\mathbf{z},\mathbf{k})$.
The term corresponding to a specific latent path is the product of
those factors that appear in the numerator of the latent path probability
and correspond to ERRW-like transitions. For every pair of states
$x,y\in\mathcal{S}\sqcup\mathbb{N}$, there are $k_{xy}$ factors,
but their sequential reinforcement depends on the order in which
$k_{xy}$ ERRW-like and $( n_{xy}-k_{xy})$ mediated transitions appear
in a specific latent path. Summing these factors over all possible
orders, one pair of states at a time, we can factorize $F(\mathbf
{z},\mathbf{k})$,
\[
F(\mathbf{z},\mathbf{k}) = \prod_{x,y:n_{xy}>0}
2^{k_{xx}\mathbf{1}_{x=y}} f_{e_{xy},\beta}(n_{xy}-\mathbf
{1}_{g(x,y)=0},k_{xy}),
\]
where
\[
f_{e,\beta} (n,k) = \sum_{u\in\{0,1\}^n, \|u\|_1=k} \prod
_{j=1}^{n} \biggl( e+(1-\beta) (j-1)+\beta\sum
_{\ell<j} u_\ell\biggr)^{u_j}
\]
and
\[
e_{xy}:= \cases{ g(x,y), &\quad if $g(x,y)>0$,
\cr
1-\beta, &\quad if
$g(x,y)=0$.}
\]

\begin{prop}
\label{recursion}
The function $f_{e,\beta}$ satisfies the following recursion for all
$0< k<n$,
%
\begin{eqnarray}
\label{soleq} f_{e,\beta}(n,k) &=& f_{e,\beta}(n-1,k) \nonumber\\[-8pt]\\[-8pt]
&&{}+f_{e,\beta}(n-1,k-1) \bigl[e-1+\beta k+(1-\beta)n\bigr],\nonumber
\end{eqnarray}
where we set, for all $n\geq0$,
\[
f_{e,\beta}(n,0) = 1 \quad\mbox{and}\quad f_{e,\beta}(n,n) =
(e)_{n\uparrow1}.
\]
\end{prop}

The recursive representation allows one to compute $p(\mathbf
{z},\mathbf{k})$ quickly. In order to obtain the values of $f_{e,\beta
}(n,k)$ for every
$n<\tilde n$, where $\tilde n$ is an arbitrarily selected integer and
$k<n$, it is sufficient to solve (\ref{soleq}) fewer than $\tilde n
^2$ times.

In the next proposition, we provide a closed-form solution for
$f_{e,\beta}$ in terms of the generalized Lah numbers, a well-known
triangular array~\cite{Comtet1974}.

\begin{defn}\label{generalizedlah}
Let $(t)_{n,V_0}$ be the generalized factorial of $t$ of order $n$ and
increments $V_0 = (v_{j})_{j\geq0}$, namely
\[
(t)_{n,V_0} = (t-v_0) (t-v_1)\cdots(t-v_{n-1})
\]
with $(t)_{0,V_0}:=1$. The \textit{generalized Lah numbers}
$C(n,k,V_0,W_0)$ (sometimes referred to as generalized Stirling
numbers), are defined by
%
\begin{equation}
\label{eqgenlah}\quad C(n,k,V_0,W_0) =\cases{1, &\quad if
$k=n=0$,
\vspace*{2pt}\cr
(w_0)_{n,V_0}, &\quad if $n>0$ and $k=0$,
\vspace*{2pt}\cr
0, &\quad
if $k>n$,
\vspace*{2pt}\cr
\displaystyle \sum_{j=0}^k
\frac{ (w_j)_{n,V_0} }{ (w_j)_{j,W_0}
(w_j)_{k-j,W_{j+1}} }, &\quad if $0< k\le n$,}\hspace*{-28pt}
\end{equation}
where $W_i = (w_{j})_{j\geq i}$.
\end{defn}

\begin{prop}\label{recursionsolution}
For any $n\geq1$ and $0<k<n$ the function $f_{e,\beta}$
coincides with
\[
f_{e,\beta}(n,k) = (1-\beta)^k C\bigl(n,n-k,V_0^{(e,\beta)},W_0^{(\beta)}
\bigr),
\]
where
%
\begin{equation}
\label{eqseqv}
V^{(e,\beta)}_0=(v_{j})_{j\geq0}:= \biggl(-\frac{e+j}{1-\beta
} \biggr)_{j\geq0}
\end{equation}
and
%
\begin{equation}
\label{eqseqw}
W^{(\beta)}_0=(w_{j})_{j\geq0}:= \biggl(j-\frac{j}{1-\beta} \biggr)_{j\geq0 }.
\end{equation}
\end{prop}


\section{Posterior simulations}
\label{bayesianinferencesection}
In this section we introduce a Gibbs algorithm for performing Bayesian
inference with the $(\theta,\alpha,\beta)$ scheme given the
trajectory of a reversible Markov chain $X_1,\ldots,X_n$.
On the basis of the almost conjugate structure of the prior model
described in the previous sections we only need to sample the latent
variables $\mathbf{k}$ conditionally on the data.
Recall that the latent variables $\mathbf{k}$ express what fraction of
the transitions in $X_1,\ldots,X_n$ are ERRW-like transitions; cf.
Figure~\ref{mechanism}.

We want to sample from $ p(\mathbf{k} \mid X_1,\ldots,X_n )$ or
equivalently from
\[
p( \mathbf{k} | \mathbf{z}) \propto p( \mathbf{k}, \mathbf{z}).
\]
Recall that $\ell$ and $\ell_x$ are functions of $( \mathbf{k},
\mathbf{z}) $ and that $\mathbf{z}=d(X_1,\ldots,X_n)$. For
simplicity, and without loss of generality, we consider the case where
initially $g(X_1,\zeta)$ is infinitesimal,\vadjust{\goodbreak} $g(\zeta,\zeta)=\theta$,
and $g(\cdot,\cdot)=0$ otherwise. The count $\ell'$ is a function of
$\mathbf{z}$ and therefore
\begin{eqnarray*}
p( \mathbf{k} | \mathbf{z}) &\propto& F(\mathbf{z},\mathbf{k})
\frac
{\prod_{x:n_x>0} (\beta(1-\alpha)^{\mathbf{1}_{x\neq X_1}}
)_{\ell_x-1\uparrow\beta} } {(\theta+\beta)_{\ell\uparrow
2\beta}}
\\
& = & F(\mathbf{z},\mathbf{k}) \frac{\prod_{x:n_x>0} (\beta(1-\alpha
)^{\mathbf{1}_{x\neq
X_1}} )_{\ell_x-1\uparrow\beta} } {
(\theta+ \ell'\beta)_{2\ell-\ell' \uparrow\beta}} \frac{(\theta+ \ell
'\beta)_{2\ell-\ell' \uparrow\beta}} {
(\theta+\beta)_{\ell\uparrow2\beta}}
\\
& \propto & F(\mathbf{z},\mathbf{k}) \frac{\prod_{x:n_x>0} (\beta(1-\alpha
)^{\mathbf{1}_{x\neq
X_1}} )_{\ell_x-1\uparrow\beta} } {
(\theta+ \ell'\beta)_{2\ell-\ell' \uparrow\beta}}
\frac{(\theta)_{2\ell\uparrow\beta}}{(\theta+\beta)_{\ell
\uparrow2\beta}}
\\
&\propto& F(\mathbf{z},\mathbf{k}) \frac{\prod_{x:n_x>0}
(\beta(1-\alpha)^{\mathbf{1}_{x\neq X_1}} )_{\ell_x-1\uparrow
\beta} } {(\theta+ \ell'\beta+\beta)_{2\ell-\ell'-1 \uparrow
\beta}} (
\theta)_{\ell\uparrow2\beta}.
\end{eqnarray*}
If
\[
G \sim\operatorname{Gamma}\bigl(\mbox{scale}=1,\mbox{shape}=\theta
/(2\beta)
\bigr)
\]
and
\[
D:= (D_{X_1}, D_1, D_2, \ldots,
D_{\ell'+1}) \sim\operatorname{Dirichlet} \biggl(1, 1-\alpha, 1-\alpha,\ldots, 1-\alpha, \frac{\theta}{\beta}+\ell'\alpha\biggr),
\]
then we can write $p(\mathbf{k} | \mathbf{z}) \propto\mathbb
{E}_{G,D} [ \psi(\mathbf{k},G,D) ]$, where
\begin{eqnarray*}
&&\psi(\mathbf{k},G,D) 
\\
&&\qquad= 
\biggl( \prod_x\frac{1}{D_x}
\biggr) \biggl(\prod_{x,y} 2^{k_{xy}\mathbf{1}_{x=y}}f_{1-\beta,\beta
}(n_{xy}-1,k_{xy})
(D_xD_y2\beta G)^{n_{xy}-k_{xy}} \biggr).
\end{eqnarray*}
In other words, if we consider the joint distribution of three
variables, $\mathbf{k}^* = \{k_{xy}^*; x,y\in\mathcal{S}\sqcup
\mathbb{N}, n_{xy}>0\}$, $D^*=(D^*_{X_1},D^*_1,D^*_2,\ldots,D^*_{\ell'+1})$
and $G^*$,
\[
p\bigl(\mathbf{k}^*, D^*,G^*\bigr) \propto p_{G}\bigl(G^*
\bigr)p_{D}\bigl(D^*\bigr) \psi\bigl(\mathbf{k}^*,G^*,D^*\bigr),
\]
where $p_{G}$ and $p_{D}$ are the distributions of $G$ and $D$, then
the marginal law of $\mathbf{k}^*$ coincides with $p(\mathbf{k} |
\mathbf{z})$. We note that sampling from $p(\mathbf{k}^*\mid
D^*,G^*)$ is simple, because the variables $k^*_{xy}$ are conditionally
independent, and that sampling from $p(D^*,G^*|\mathbf{k}^*)$ is
straightforward. The random variables $D^*$ and $G^*$ conditionally on
$\mathbf{k}^*$ are independent with Dirichlet and Gamma distributions.

Finally, we use these conditional distributions to construct a Gibbs
sampler for $p(\mathbf{k}^*,D^*,G^*)$. In any Markov chain Monte Carlo
algorithm, it is important to ensure mixing. In Appendix D, we derive
an exact sampler for $p(\mathbf{k}|\mathbf{z})$ which uses a coupling
of the Gibbs Markov chain just defined. The method is related to
Coupling From The Past~\cite{Propp1996}. We performed simulations
with the exact sampler to check the convergence of the proposed Gibbs algorithm.


\section{Analysis of molecular dynamics simulations}
\label{applications}

\subsection{The species sampling problem}

Species sampling problems have a long history in ecological and
biological studies. The aim is to determine the species composition of
a population containing an unknown number of species when only a sample
drawn from it is available.

A common statistical issue is how to estimate species richness, which
can be quantified in different ways. For example, given an initial
sample of size $n$, species richness might be quantified by the number
of new species we expect to observe in an additional sample of size
$m$. It can be alternatively evaluated in terms of the probability of
discovering at the $(n+m)$th draw a new species that does not appear
across the previous $(n+m-1)$ observations; this yields the discovery
rate as a function of the size of an hypothetical additional sample.
These estimates allow one to infer the coverage of a sample of size
$n+m$, in other words, the relative abundance of distinct species
observed in a sample of size $n+m$.

A review of the literature on this problem can be found in Bunge and
Fitzpatrick~\cite{Bunge1993}. Lijoi et al. proposed a Bayesian
nonparametric approach for evaluating species richness, considering a
large class of exchangeable models, which include as special case the
two-parameter Hoppe urn~\cite{Lijoi2007a}. See also Lijoi et al.
\cite{Lijoi2007b} and Favaro et al.~\cite{Favaro2009} for a
practitioner-oriented illustration using expressed sequence tag (EST)
data obtained by sequencing cDNA libraries.

We illustrate the use of the $(\theta,\alpha,\beta)$ scheme in
species sampling problems. In particular, we evaluate species richness
in molecular dynamics simulations.

\subsection{Data}
The data we analyze come from a series of recent studies applying
Markov models to protein molecular dynamics simulations (\cite{Pande2010}, and
references therein). These computer experiments produce
time series of protein structures. The space of structures is
discretized, such that two structures in a given state are
geometrically similar; this yields a sequence of species which
correspond to conformational states that the molecule adopts in water.
We apply the $(\theta,\alpha,\beta)$ scheme to perform predictive
inference of this discrete time series.

We analyze two datasets. The first is a simulation of the alanine
dipeptide, a~very simple molecule. The dataset consists of 25,000
transitions, sampled every 2 picoseconds, in which 104 distinct states
are observed. In this case the 50 most frequently observed states
constitute 85\% of the chain and each of the 104 observed states
appears at least 12 times. The second dataset is a simulation of a more
complex protein, the WW domain, performed in the supercomputer Anton
\cite{Shaw2010}. This example illustrates the complexity of the
technology and the large amount of resources required for simulating
protein dynamics \textit{in silico}. It also motivates the need for
suitable statistical tools for the design and analysis of these
experiments. In this dataset 1410 distinct states are observed in
10,000 transitions, sampled every 20 nanoseconds. Many of the states
are observed only a few times; in particular we have 991 states that
have been observed fewer than 4 times and 547 states that appear only once.

\subsection{Prior specification}

To apply the $(\theta,\alpha,\beta)$ scheme it is necessary to tune
the three parameters. We consider the initial weights $g$ everywhere
null except for $g(\zeta,\zeta)=\theta$ and $g(X_1,\zeta)$
infinitesimal. The parameters $\theta$ and $\alpha$ affect the
probability of finding a novel state when the latent process reaches
$\zeta$, while the parameter $\beta$ tunes the degree of dependence
between the random transition probabilities. We recall that in the
extreme case of $\beta=1$ the sequence $(X_{i})_{i\in\mathbb{N}}$ is
exchangeable and the random transition probabilities out of the
observed states become identical.

We proceed by approximating the marginal likelihood of the data for the
set of parameters $(\theta,\alpha,\beta)$, where
$\theta\in\{1,5,10,25,50,100,300,400,500\}$,
$\alpha\in\{0.03,0.2,0.5,0.8,0.97\}$ and
$\beta\in\{0.03,0.2,0.5,0.8,0.97\}$.
We iteratively drew samples, under specific $(\theta,\alpha,\beta)$
values, from the conditional distribution $p(\mathbf{k} | \mathbf
{z})$ using the Gibbs algorithm defined in the previous section. Note that
\[
\eta:= \sum_{\mathbf{k}} \frac{1}{p(\mathbf{z},\mathbf{k})} p(\mathbf
{k} |
\mathbf{z}) = \sum_{\mathbf{k}} \frac{1}{p(\mathbf
{z},\mathbf{k})}
\frac{p(\mathbf{z},\mathbf{k})}{p(\mathbf{z})} = \frac{\prod_{x,y\in
\mathbb{N}\sqcup\mathcal{S}} (n_{xy}-1)
}{p(\mathbf{z})},
\]
where the last equality is obtained counting the possible values of
$\mathbf{k}$.
We compare the models on the basis of approximations $\hat p(\mathbf
{z})$ of the marginal probabilities $p(\mathbf{z})$ across prior
parameterizations. The samples we drew from $p(\mathbf{k} | \mathbf
{z})$ are used to compute Monte Carlo estimates $\hat{\eta}$ of $\eta
$. Recall that the probability $p(\mathbf{k}, \mathbf{z})$ can be
computed using the analytic expressions derived in Section \ref
{lawsection}. Using $\hat\eta$ we compute the estimates $\hat{p}(\mathbf
{z}):= \prod_{x,y\in\mathbb{N}\sqcup\mathcal{S}} (n_{xy}-1) /
\hat{\eta}$ and obtain standard errors by Bootstrapping.

In Table~\ref{modelcomparisons}, we report the logarithm of these
estimates for each model, shifted by a constant such that the largest
entry for each dataset is 0. We only show, due to limits of space,
these results for the $\theta$ values associated with the maxima of
$\hat p(\mathbf{z})$ across the considered parameterizations. The
difference between two entries corresponds to a logarithmic Bayes
factor between two models. The values in Table~\ref{modelcomparisons}
indicate that in each dataset there is one model for which there is
strong evidence against all others. This also holds when several values
of $\theta$ are considered. For each dataset, we have highlighted the
optimal parameters. The degenerate cases $\alpha=0$ and $\beta=1$
were also included in the comparisons but are not shown in Table \ref
{modelcomparisons}. The difference in the marginal log-likelihood
between models with $\alpha=0$ and $\alpha=0.03$ is negligible. On
the other hand, shifting the parameter $\beta$ from 0.97 to 1 in the
optimal model for dataset 2 decreased the log-likelihood by 7565, as
this model is exchangeable and does not capture the Markovian nature of
the data. These observations suggest that a fully Bayesian treatment
with a hyper-prior over a grid of possible $(\theta,\alpha,\beta)$
combinations would produce similar results.

\begin{table}
\caption{The log-likelihood of the data, $\log p(\mathbf{z})$, for
the $(\theta,\alpha,\beta)$ scheme at the optimal value of $\theta$
for a range of values of $\alpha$ and $\beta$. In each table, every
entry is shifted by a constant such that the largest entry equals 0}
\label{modelcomparisons}
\begin{tabular*}{\tablewidth}{@{\extracolsep{\fill}}llrcrl@{}}
\hline
& \multicolumn{5}{c@{}}{$\bolds{\beta}$}\\[-4pt]
& \multicolumn{5}{c@{}}{\hrulefill}\\
$\bolds{\alpha}$ & \multicolumn{1}{c}{\textbf{0.03}}
& \multicolumn{1}{c}{\textbf{0.2}} & \multicolumn{1}{c}{\textbf{0.5}}
& \multicolumn{1}{c}{\textbf{0.8}} & \multicolumn{1}{c@{}}{\textbf{0.97}}\\
\hline
\multicolumn{6}{@{}c@{}}{Dataset 1: Alanine dipeptide
($\theta=25$)} \\
[4pt]
0.03 & $-3212\pm 10.4$ & \multicolumn{1}{l}{\hphantom{000}$-170\pm 11.7\hspace*{-4.5pt}$}
& \hphantom{0$-$}\fbox{$0\pm 5.9$} &
$-982\pm 7.4$ & $-2828\pm 12.3$ \\
0.2 & $-3263\pm 4.0$ & $-220\pm 6.6$ & $-28\pm 8.6$ & $-997\pm 6.0$
& $-2809\pm 9.5$ \\
0.5 & $-3404\pm 2.1$ & $-333\pm 4.1$ & $-125\pm 15.5$ & $-1024\pm 5.4$
& $-2815\pm 2.3$ \\
0.8 & $-3621\pm 6.9$ & $-525\pm 6.9$ & $-232\pm 3.0$\hphantom{0} & $-1099\pm 4.2$
& $-2857\pm 5.8$ \\
0.97 & $-3868\pm 3.1$ & $-763\pm 4.6$ & $-447\pm 3.4$\hphantom{0} & $-1280\pm 4.2$
& $-2960\pm 11.9$ \\
[4pt]
\multicolumn{6}{@{}c@{}}{Dataset 2: WW domain ($\theta=500$)} \\
[4pt]
0.03 & $-14\mbox{,}695\pm 2.0$ & $-5147\pm 2.2$ & $-1701\pm 1.8$ & $-361\pm 4.0$
& $-234\pm 1.9$ \\
0.2 & $-15\mbox{,}167\pm 1.6$ & $-5507\pm 1.4$ & $-1865\pm 4.2$ & $-329\pm 3.5$
& \hphantom{0}$-95\pm 3.2$ \\
0.5 & $-16\mbox{,}211\pm 2.4$ & $-6354\pm 2.4$ & $-2365\pm 0.8$ & $-482\pm 5.6$
& \hphantom{0$-$\hspace*{1.2pt}}\fbox{$0\pm 1.2$} \\
0.8 & $-17\mbox{,}943\pm 1.6$ & $-7893\pm 1.8$ & $-3542\pm 3.0$ & $-1120\pm 7.6$
& $-119\pm 1.4$ \\
0.97 & $-20\mbox{,}892\pm 1.9$ & $-10\mbox{,}739\pm 5.0$ & $-6143\pm 1.5$ & $-3194\pm 7.1$
& $-964\pm 0.9$\\
\hline
\end{tabular*}
\end{table}

Summarizing, the use of a three-dimensional grid and the computation of
Monte Carlo estimates allows one to effectively obtain a parsimonious
approximation of the likelihood function that, in our case, supported
selection of single parameterizations.

\subsection{Posterior estimates}

\begin{figure}

\includegraphics{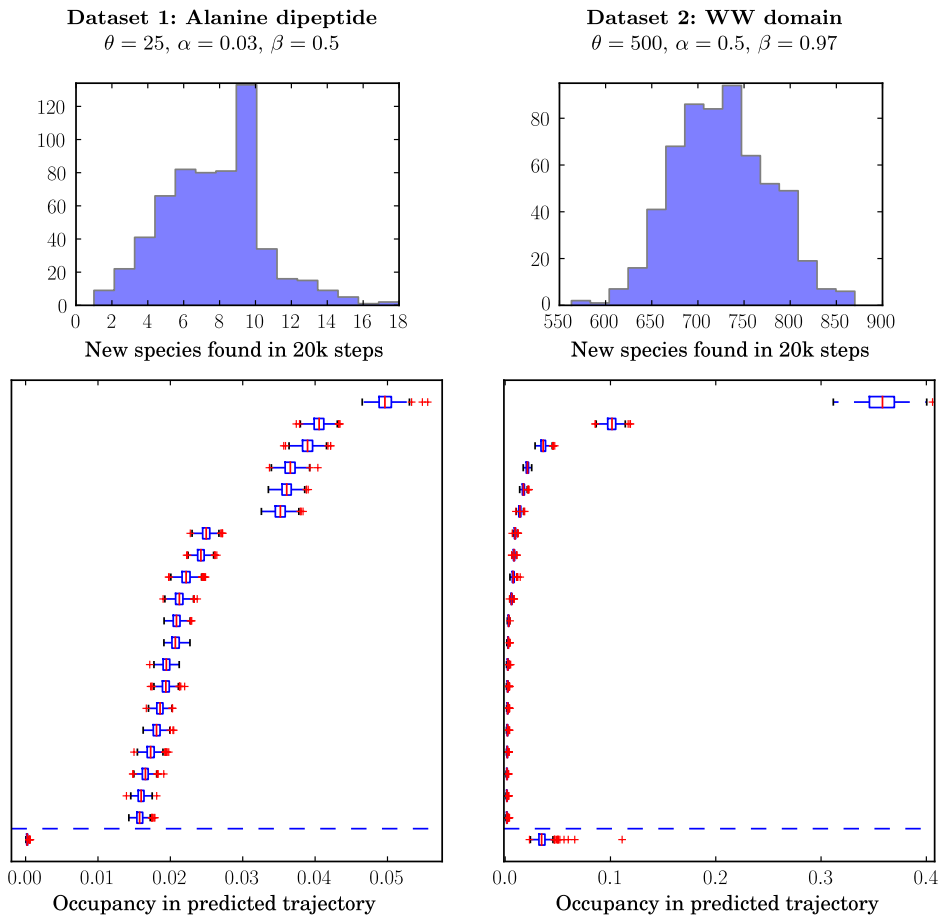}

\caption{Posterior simulations for two molecular dynamics datasets.
Top: histogram of the number of new species found in 600 simulations
from the predictive distribution for $X_{n+1},\ldots,X_{n+20\mbox{,}000}$,
given the data $X_1,\ldots,X_n$. Bottom: box plot of the fraction of
time spent at each state in these simulations. Only the twenty most
populated states are shown; below the dashed line, we show the fraction
of time spent at states not observed in the dataset.}
\label{posteriorplots}
\end{figure}

The main results of our analysis are summarized in Figure \ref
{posteriorplots}. Conditional on each sample $\mathbf{k} \sim
p(\mathbf{k}|\mathbf{z})$, generated under the selected
$(\theta,\alpha,\beta)$ parametrization, we simulated 20,000 future
transitions using our predictive scheme. Once $\mathbf{k}$ is
%
conditionally sampled, the predictive simulations become
straightforward with the reinforcement scheme. To provide a measure of
species richness and the associated uncertainty, we histogram the
number of new states discovered in our simulations in Figure \ref
{posteriorplots}. Only a few states are predicted to be found for
dataset 1, while a large number of new states are predicted for
dataset~2. This result is not surprising because the alanine dipeptide
dataset has a limited number of rarely observed states, while in the WW
domain data a significant number of states are observed once. This
result also seems consistent with the selected values of $\theta$ and
$\beta$ in these two experiments.

As previously mentioned the $(\theta,\alpha,\beta)$ scheme is a
Bayesian tool for predicting any characteristic of the future
trajectories $X_{n+1},\ldots, X_{n+m}$. The bottom panels in Figure
\ref{posteriorplots} show confidence bands for the predicted
fractions of time that will be spent at the most frequently observed
states in the next 20,000 transitions. Each box in the plots refers to
a single state and shows the quartiles and the 10th and 90th
percentiles of the predictive distribution; states are ordered
according to their mean observed frequency. We only show these
occupancies for the 20 most populated states, and below the dashed
line, we show the total occupancy for states that do not appear in the
original data. In the WW domain example, the simulation is expected to
spend between 2.5\% and 5\% of the time at new states.

\begin{figure}

\includegraphics{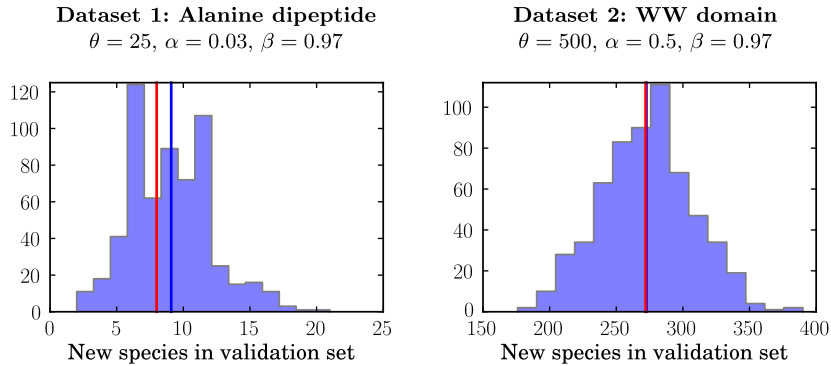}

\caption{Bayesian predictions for the number of new species in a
validation set. The histograms show the number of species found in 600
simulations from the predictive distribution for
$X_{n+1},\ldots,X_{n+m}$, where $m$ is the length of the validation set
and $X_1,\ldots,X_n$ is the training set. The blue line shows the mean
of these samples. The red line shows the actual number of new species
found in the validation set. Note that in the right panel, the lines
overlap due to the small separation between them.}
\label{posteriorcrossvalidation}
\end{figure}

To assess the predictive performance of the model we split each dataset
into a training set and a validation set. The rationale of this
procedure is identical to routinely performed cross validations for
i.i.d. data. In our setting, the training and validation sets are
independent portions of a homogeneous Markov chain. The first part of
the procedure, which uses only the training set, includes selection of
the $(\theta,\alpha,\beta)$ parameters and posterior computations.
Then, we contrast Bayesian predictions to statistics of the validation
set. Overall, this approach suggests that our model generates reliable
predictions. Figure~\ref{posteriorcrossvalidation} shows histograms
for the number of new species found in predictive simulations of equal
length as the validation set. In each panel, the blue line is the Bayes
estimate and the red line corresponds to the number of species that was
actually discovered in the validation set. This approach also supports
the inference reported with box plots in Figure~\ref{posteriorplots}.
We repeated the computations for deriving the results in Figure \ref
{posteriorplots} using only the training data, and considering a
future trajectory equal in length to the validation data.
In this case, 37 out of 42 of the true state occupancies in the
validation set were contained in the 90\% posterior confidence bands.


\section{Discussion}
\label{discussionsection}

We introduced a reinforced random walk with a simple predictive
structure that can be represented as a mixture of reversible Markov
chains. The model generalizes exchangeable and partially exchangeable
sequences that have been extensively studied in the literature. Our
nonparametric prior, the de Finetti measure of the
$(\theta,\alpha,\beta)$ scheme, can be viewed as a distribution over
weighted graphs with a countable number of vertices in a possibly
uncountable space $\mathcal{X}$. As is the case for other well-known
Bayesian nonparametric models such as the Dirichlet process~\cite
{Ferguson1973}, the hierarchical Dirichlet process~\cite{Teh2006} and
the infinite hidden Markov model~\cite{Beal2002}, it is possible to
represent our model as a function of two independent components, a
species sampling sequence $(Z_i)_{i \in\mathbb{N}}$ and a process
$(T_i)_{i \in\mathbb{N}}$ which determines the species' locations. This
property is fundamental in applications including Dirichlet process
mixture models and the infinite hidden Markov model.

A natural extension of our model, not tackled here, is the definition
of hidden reversible Markov models. A simple construction would consist
of convolving our vertices with suitable density functions. We hope
reversibility can be an advantageous assumption in relevant
applications; in particular we think reversibility can be explored as a
tool for the analysis of genomic data and time series from
single-molecule biophysics experiments.

\begin{appendix}\label{app}
\section*{\texorpdfstring{Appendix: Proofs from Sections \protect\lowercase{\ref{tabschemesection}} and \protect\lowercase{\ref{representationssection}}}
{Appendix: Proofs from Sections 2 and 3}}

\begin{pf*}{Proof of Proposition~\ref{recurrence}}
Consider the latent process on $\mathcal{X}_+$.
The transition probability from $x_1\in\mathcal{X}_+$ to $x_2\in
\mathcal{X}_+$ with $g(x_1,x_2)>0$ is of the form $g(x_1,x_2)/\sum_{y\in
\mathcal{X}_+} g(x_1,y)$.
Between successive visits to $x_1$, the denominator $\sum_{y\in
\mathcal{X}_+} g(x_1,y)$ is increased by at most $2$, and the
numerator may only increase.
Assume that almost surely, the process visits $x_1$ infinitely often.
There exist $c_2\ge c_1>0$, such that if $A_{n,m}$ is the event that we
do not traverse $(x_1,x_2)$ between the $n$th and $m$th visits to $x_1$,
\[
p(A_{n,m}) \leq\prod_{i=n}^{m-1} \biggl[ 1- \frac
{c_1}{c_2+2(i-n)} \biggr],
\]
which goes to 0 as $m\to\infty$.
Therefore the edge $(x_1,x_2)$ is a.s. traversed infinitely often.
Thus, if $x_1\in\mathcal{X}_{+}$ is a.s. visited infinitely often, by
induction the process a.s. returns infinitely often to all visited
states. Suppose a state in $\mathcal{X}$ is visited infinitely often
a.s., then the process visits $\zeta$ infinitely often by the previous
argument. Otherwise, the process must visit an infinite number of
states in $\mathcal{X}$, and since the set of pairs $(x,y)\in\mathcal
{X}^2$ with a positive initial weight $g(x,y)$ is a finite subset of
$\mathcal{S}^2$, we must go through $\zeta$ an infinite number of times.
We conclude that $\zeta$ is visited infinitely often a.s. and
therefore the process returns to every state visited infinitely often.
If $\theta>0$, then the edge $(\zeta,\zeta)$ is crossed infinitely
often, and we see an infinite number of distinct states.
\end{pf*}


\begin{pf*}{Proof of Proposition~\ref{partialexchangeability}}
For $1\leq i< n$ there is a latent variable $U_i\in\{a,b,c\}$ that
determines in which of the three ways outlined in Figure
\ref{mechanism} the transition\vadjust{\goodbreak} $X_i \to X_{i+1}$ proceeded. The
probability of $Z_1=x_1,\ldots,Z_n=x_n$ is the sum of its joint
probability with every latent sequence $U_1=u_1,\ldots,U_{n-1}=u_{n-1}$.
We will show that there is a one-to-one map $L$ of the latent sequences
such that, letting $L(u_1,\ldots,u_{n-1})=u'_1,\ldots,u'_{n-1}$, one has
%
\begin{eqnarray}
\label{ad}
&&
p(Z_{1}=x_1,\ldots,Z_{n}=x_n;U_{1}=u_1,\ldots,U_{n-1}=u_{n-1})
\nonumber\\[-8pt]\\[-8pt]
&&\qquad =
p\bigl(Z_{1}=x'_1,\ldots,Z_{n}=x'_n;U_{1}=u_{1}',\ldots,U_{n-1}=u_{n-1}'\bigr).
\nonumber
\end{eqnarray}
The proposition follows from this claim. Let $t$ denote the
transposition such that $d(t(x_1,\ldots,x_n))=x'_1,\ldots,x'_n$. Define
$x''_1,\ldots,x''_n=t(x_1,\ldots,x_n)$. The map $L$ is defined so that
for any $u_{1},\ldots,u_{n-1}$ and $j,k<n$, satisfying
\begin{eqnarray*}
\mathbf{1}_{(x_j=x, x_{j+1}=y)}+ \mathbf{1}_{(x_j=y, x_{j+1}=x)}&>&0,
\\
\mathbf{1}_{(x''_k=x, x''_{k+1}=y)}+ \mathbf{1}_{(x''_k=y, x''_{k+1}=x)}&>&0
\end{eqnarray*}
for some $x$ and $y$, if
\[
\sum_{i=1}^j \mathbf{1}_{(x_i=x, x_{i+1}=y)}+
\mathbf{1}_{(x_i=y,
x_{i+1}=x)}= \sum_{i=1}^k
\mathbf{1}_{(x''_i=x, x''_{i+1}=y)}+ \mathbf{1}_{(x''_i=y,
x''_{i+1}=x)},
\]
then $u_j=u'_k$. Note that we can define the joint probability of
$(Z_i)_{i\le n}$ and $(U_i)_{i<n}$ through a reinforcement scheme
identical the one defined in Section~\ref{tabschemesection}.
Precisely, the probability of each transition and associated
category is of the form
\[
\frac{g(x_i,x_{i+1})}{g(x_i)} \mathbf{1}_{u_i=a},\qquad \frac{g(x_i,\zeta
)g(\zeta,x_{i+1})}{g(x_i)g(\zeta)}\mathbf
{1}_{u_i=b},\qquad \frac{g(x_i,\zeta)g(\zeta,\zeta)}{g(x_i)g(\zeta
)}\mathbf{1}_{u_i=c},
\]
where {$g(x):= \sum_{y\in\mathcal{X}_{+}} g(x,y)$.}

The factors $g(\zeta)$, which appear in the denominator when $U_i\in\{
b,c\}$, are reinforced by $2\beta$ between successive visits to $\zeta
$. Therefore, their product only depends on the number of mediated
transitions, which is invariant under~$L$. Similarly, factors
$g(\zeta,\zeta)$ increase by $\alpha\beta$ between successive
occurrences; their product is identical when we compute the two sides
of (\ref{ad}) because the number of mediated transitions with discovery
remains identical. Also, the factors $g(x)$ in the denominators
increase by $2$ between successive occurrences of the same
$\mathcal{X}$ state; their product is identical when we compute the two
sides of (\ref{ad}) because the number of transitions out of any state
$x$ (or toward $x$) remains identical. Finally, we need to prove the
identity between
%
\begin{equation}
\label{ads1} \prod_{i} \bigl(g(x_i,x_{i+1})
\mathbf{1}_{u_i=a} +g(x_i,\zeta)g(x_{i+1},\zeta)
\mathbf{1}_{u_i=b}+g(x_i,\zeta)\mathbf{1}_{u_i=c}
\bigr)
\end{equation}
and
%
\begin{equation}
\label{ads2} \prod_{i} \bigl(g
\bigl(x'_i,x'_{i+1}\bigr)
\mathbf{1}_{u'_i=a} +g\bigl(x'_i,\zeta\bigr)g
\bigl(x'_{i+1},\zeta\bigr)\mathbf{1}_{u'_i=b}+g
\bigl(x'_i,\zeta\bigr)\mathbf{1}_{u'_i=c}
\bigr).
\end{equation}
The identity between (\ref{ads1}) and (\ref{ads2}) follows by
combining the definitions of $d$ and $L$ with the reinforcement
mechanism. Specifically, the factors $g(x,\zeta)$ and $g(x',\zeta)$
are increased by $\beta$ between successive occurrences. Since
$g(x,\zeta)$ appears as many times in the left-hand side of (\ref
{ad}) as $g(x',\zeta)$ does in the right-hand side of (\ref{ad}), the
product of these factors is identical in each case. The remaining
factors $g(x,y)$ may increase by different amounts between successive
occurrences. Their product is a function of the subsequence of
$U_{1},\ldots,U_{n-1}$ with indices $\{1\leq i< n\dvtx  \{Z_{i},Z_{i+1}\}=\{
x,y\}\}$. By the definition of $L$, this subsequence is the same in the
left and right-hand sides of (\ref{ad}), which completes the proof of
our claim.
\end{pf*}

\begin{pf*}{Proof of Proposition~\ref{representationtheorem}}
Let $(X_i^{\prime\prime})_{i\in\mathbb{N}}$ be a
$(\theta,\alpha,\beta)$ scheme. The process
$(X_i^{\prime\prime})_{i\in\mathbb{N} }$ returns to
$X_1^{\prime\prime}$ infinitely often a.s. Let $h_i$ be the $i$th
$X_1^{\prime\prime}$-block. Define $(X_i)_{i\in\mathbb{N}
}:=(h_1,h_3,h_5,\ldots) $ and $(X^{ \prime}_i)_{i\in\mathbb{N}
}:=(h_2,h_4,h_6,\ldots)$.~Proposi-\break tion~\ref{partialexchangeability}
implies $(X_i)_{i\in\mathbb {N} }\stackrel{\mathrm{d}}{=}
(X^\prime_i)_{i\in\mathbb{N} }\stackrel{\mathrm{d}}{=}
(X^{\prime\prime}_i)_{i\in\mathbb{N} }$. Let $\operatorname{Fin}(x)$ be the
last element of a vector $x$. Define
\[
W_i:= \lim_{m\rightarrow\infty} \operatorname{Fin}\bigl( d
\bigl(X^{\prime
}_1,X^{\prime}_2,\ldots,X^{\prime}_m,X_1,X_2,\ldots,X_i\bigr) \bigr).
\]
This limit exists a.s. because $X_i$ is recurrent in $(X''_j)_{j\in
\mathbb{N}}$; therefore, the sequence of blocks that form
$(X'_j)_{j\in\mathbb{N}}$ is conditionally i.i.d. from a distribution
which a.s. assigns positive probability to blocks containing $X_i$,
which implies $(X'_j)_{j\in\mathbb{N}}$ visits $X_i$ after a finite
time a.s., at which point the limit settles.

In Lemma~\ref{markovexchangeabilityofW} we show that $(W_i)_{i\in
\mathbb{N}}$ is Markov exchangeable and recurrent. Therefore, by de
Finetti's theorem for Markov chains (\ref
{definettistheoremformarkovchains}), it is a mixture of Markov chains.
Finally, by Lemma
\ref{iidorderdoesntmatter}, we obtain the representation claimed
in the proposition.
\end{pf*}

\begin{lem}
\label{markovexchangeabilityofW}
Without loss of generality, let $\mathcal{X}=(0,1)$. The process
$(W_i)_{i\in\mathbb{N}}$ is Markov exchangeable and returns to every
state in $\mathcal{S}\sqcup\mathbb{N}$ infinitely often a.s.
\end{lem}

\begin{pf}
The recurrence of $(Z_i)_{i\in\mathbb{N}}$, which is a consequence of
Proposition~\ref{recurrence}, implies the recurrence of $(W_i)_{i\in
\mathbb{N}}$. Thus, we have left to show Markov exchangeability.

The sequence $W_1,\ldots,W_n$ can be mapped through $d$ to
$Z_1,\ldots,Z_n$, which is a species sampling sequence for the
$(\theta,\alpha,\beta)$ scheme. Take any sequence $w_1,\ldots,w_n$ and
let $z_1,\ldots,z_n:= d(w_1,\ldots,w_n)$. We have
%
\begin{eqnarray}\label{probw}
&&
p(W_1=w_1,\ldots,W_n=w_n)\nonumber\\
&&\qquad=
p(Z_1=z_1,\ldots,Z_n=z_n)
\\
&&\qquad\quad{}\times p(W_1=w_1,\ldots,W_n=w_n
| Z_1=z_1,\ldots,Z_n=z_n).\nonumber
\end{eqnarray}

Consider any pair of sequences $w_1,\ldots,w_n$ and $w'_1,\ldots,w'_n$
related by a transposition of two blocks with identical\vadjust{\goodbreak} initial and
final states.
Proposition~\ref{partialexchangeability} implies
\[
p \bigl((Z_1,\ldots,Z_n)=d(w_1,\ldots,w_n) \bigr)=p
\bigl((Z_1,\ldots,Z_n)=d \bigl(w'_1,\ldots,w'_n\bigr) \bigr).
\]
We have left to show that the second factor on the right-hand side of
(\ref{probw}) is identical for $w_1,\ldots,w_n$ and
$w'_1,\ldots,w'_n$.
The identity of the conditional distribution of $(X_i^{\prime})_{i\in
\mathbb{N} }$ given $(Z_1,\ldots,Z_n)$ equal to $d(w_1,\ldots,w_n)$
or equal to $ d(w'_1,\ldots,w'_n) $
proves the lemma.
\end{pf}

\begin{lem}
\label{iidorderdoesntmatter}
The process $(T_{W_i}\mathbf{1}_{W_i\notin\mathcal{S}} + W_i\mathbf
{1}_{W_i\in\mathcal{S}})_{i\in\mathbb{N}}$ has the same
distribution as $(X_i)_{i\in\mathbb{N}}$.
\end{lem}

\begin{pf}
By definition $(X_i)_{i\in\mathbb{N}}\stackrel{\mathrm{d}}{=}
(T_{Z_i}\mathbf{1}_{Z_i\notin\mathcal{S}} + Z_i\mathbf{1}_{Z_i\in
\mathcal{S}})_{i\in\mathbb{N}}$. Note that\break $(T_i)_{i\in\mathbb
{N}}$ is an i.i.d. sequence, independent from $(Z_i)_{i\in\mathbb
{N}}$, and $d((W_i)_{i\in\mathbb{N}}) = (Z_i)_{i\in\mathbb{N}}$.
These facts imply that $(T_{W_i}\mathbf{1}_{W_i\notin\mathcal{S}} +
W_i\mathbf{1}_{W_i\in\mathcal{S}})_{i\in\mathbb{N}}
\stackrel{\mathrm{d}}{=} (T_{Z_i}\mathbf{1}_{Z_i\notin\mathcal{S}} +
Z_i\mathbf{1}_{Z_i\in\mathcal{S}})_{i\in\mathbb{N}}$.
\end{pf}

\begin{pf*}{Proof of Proposition~\ref{reversibility}}
Let $h_1,h_2,\ldots$ be the $X_1$-blocks of the $(\theta,\alpha,\beta)$
scheme. Consider a map $s$ on the $X_1$-blocks' space; if
$a=(a_1,\ldots,a_m)$, then $s(a)=(a_1,a_m,a_{m-1},\ldots, a_2)$. We
can observe, following the same arguments used for proving Proposition
\ref{partialexchangeability}, that for any $X_1$-block $a$ and any
integer $n$,
\[
p(h_1=a, h_2,\ldots,h_n)=p
\bigl(h_1=s(a),h_2,\ldots,h_n\bigr).
\]

Let $F$ be a random measure distributed according to the de Finetti
measure of the $X_1$-blocks.
The above expression and the equality
\[
\lim_{n\to\infty} p(h_1\in A\mid h_2,\ldots,h_n)\stackrel{\mathrm{a.s.}} {=}F(A),
\]
where $A$ is a generic measurable set,
imply that a.s. the distance in total variation between $F$ and $F
\circ s$
is null.
\end{pf*}
\end{appendix}

\section*{Acknowledgments}

We are grateful to an Associate Editor and three Referees for their
constructive comments and suggestions. We would like to thank DE Shaw
Research~\cite{Shaw2010} for providing the molecular dynamics
simulations of the WW domain. The Markov models analyzed in Section
\ref{applications} were generated by Kyle Beauchamp, using the
methodology described in~\cite{Pande2010}. We would also like to thank
Persi Diaconis and Vijay Pande for helpful suggestions.

\begin{supplement}
\stitle{Appendices B, C and D}
\slink[doi]{10.1214/13-AOS1102SUPP} 
\sdatatype{.pdf}
\sfilename{aos1102\_supp.pdf}
\sdescription{Appendix~B describes the two-parameter HDP-HMM in
relation to the $(\theta,\alpha,\beta)$ scheme. Appendix C contains all
proofs from Sections~\ref{largesupportsection}, \ref
{sufficientnesssection} and~\ref{lawsection}. Appendix D contains a
derivation of
the exact sampler mentioned in Section~\ref{bayesianinferencesection}
using Coupling From the Past.}
\end{supplement}


\printaddresses

\end{document}